\input amstex
\documentstyle{amsppt}

\footline{}

\catcode`\@=11

\newfam\ssffam
\font\tenssf=cmss10
\font\eightssf=cmss8
\def\ssf{\fam\ssffam}
  \textfont\ssffam=\tenssf \scriptfont\ssffam=\eightssf
\def\sf#1{\leavevmode\skip@\lastskip\unskip\/%
       \ifdim\skip@=\z@\else\hskip\skip@\fi{\ssf#1}}
\def\rom#1{\leavevmode\skip@\lastskip\unskip\/%
        \ifdim\skip@=\z@\else\hskip\skip@\fi{\rm#1}}
\def\newmcodes@{\mathcode`\'"27\mathcode`\*"2A\mathcode`\."613A%
 \mathcode`\-"2D\mathcode`\/"2F\mathcode`\:"603A }
\def\operatorname#1{\mathop{\newmcodes@\kern\z@\fam\z@#1}}
\catcode`\@=\active

\catcode`\@=11
\def\raggedcenter@{\leftskip\z@ plus.3\hsize \rightskip\leftskip
 \parfillskip\z@ \parindent\z@ \spaceskip.3333em \xspaceskip.5em
 \pretolerance9999\tolerance9999 \exhyphenpenalty\@M
 \hyphenpenalty\@M \let\\\linebreak}
\catcode`\@=\active

\mag=1200
\pagewidth{14truecm}
\pageheight{21.5truecm}
\hoffset23pt
\voffset36pt
\binoppenalty=10000
\relpenalty=10000
\tolerance=500
\mathsurround=1pt

\define\m1{^{-1}}
\define\ov1{\overline}
\def\gp#1{\langle#1\rangle}
\def\ul2#1{\underline{\underline{#1}}}
\def\Y{\mathbin{\ssf Y@!@!@!@!}}

\hrule height0pt
\vskip 1truein

\topmatter
\title 
On the order of the unitary subgroup of modular group algebra
\endtitle

\author
  Victor Bovdi\\
{\eightpoint 
Lajos Kossuth University, H-410 Debrecen, P.O. Box 12, Hungary 
}\\
{\eightpoint 
      vbovdi\@math.klte.hu 
}
\phantom{anything}\\
\phantom{anything}\\
  A.\,L.\,Rosa\\
{\eightpoint 
Universidade Federal de Ouro Preto,
35400-000 Ouro Preto--M.G., Brasil
}\\
{\eightpoint 
alrosa\@iceb.ufop.br 
}
\endauthor

\leftheadtext\nofrills{ BOVDI AND ROSA }
\rightheadtext\nofrills{
ON THE ORDER OF THE UNITARY SUBGROUP ...
}

\abstract
Let $KG$ be a group algebra of a finite $p$-group $G$ over a finite
field $K$ of characteristic $p$. We compute the order of the
unitary subgroup of the group of units  when $G$ is either  an
extraspecial $2$-group or the central product of such a group with
a cyclic group of order $4$ or $G$ has an abelian subgroup $A$ of
index $2$ and an  element $b$  such that $b$ inverts each element
of $A$.
\endabstract


\thanks 
The first author's research was supported  by the Hungarian 
National Foundation for Scientific Research No.T 025029  
and by FAPESP Brazil (proc. 97/05920-6).
\endthanks

\endtopmatter
\document

\subhead
1. Introduction
\endsubhead
Let $KG$ be the group algebra of a locally-finite $p$-group $G$
over a commutative ring $K$ (with~$1$) and $V(KG)$ be the group of
normalized units (that is, of the units with augmentation $1$) in
$KG$. The anti-automorphism $g\mapsto g\m1$ extends linearly to an
anti-automorphism $a\mapsto a^*$ of $KG$; this extension leaves
$V(KG)$ setwise invariant and its restriction to $V(KG)$ followed
by $v\mapsto v\m1$ gives an automorphism of $V(KG)$. The elements
of $V(KG)$ fixed by this automorphism are the {\it unitary
normalized units} of $KG$; they form a subgroup,  which we denote
by $V_*(KG)$. Interest in unitary units arose in algebraic topology
and a more general definition, involving an `orientation
homomorphism', is also current; the special case we use here arises
when the orientation homomorphism is trivial.

In \cite{3-4} A.\,Bovdi and A.\,Szak\'{a}cs solved the problem,
posed by S.~P.~No\-vikov  of the structure of the group $V_*(KG)$ of
the group  algebra of a  finite abelian $p$-group over a finite
field of $p^m$ elements.  We also know a few facts on $V_*(KG)$
when $G$ is   nonabelian ( see survey \cite{1}). Note that A.~Bovdi
and L.~Erdei \cite{2} have described the unitary subgroup $V_*(KG)$
for all groups of order $8$ and $16$.  We shall  study  here the
order of the group $V_*(KG)$  for some nonabelian groups $G$.

\subhead
2. Basic facts 
\endsubhead 
For an arbitrary element $x=\sum_{g\in G}\alpha_gg\in KG$ we 
put $\chi(x)=\sum_{g\in G}\alpha_g\in K$, for $a,b\in G$ we denote 
$a^b=b\m1ab$ and let $|a|$ denote the order of element $a$.

First, let us  recall  some facts about the order of $V_*(KG)$ for
$char(K)>2$, which are known  from the theory of algebras with
involution.

We  use the following notation:
$$
I_K(G)=\{\sum_{g\in G}\alpha_gg\in KG \,\, 
\mid \,\, \sum_{g\in G}\alpha_g=0\}
$$
for the  augmentation ideal of $KG$. Let 
$$
S(I)=\{x\in I_K(G) \,\, \mid \,\, x^*=x  \}, \,\,\,\,\,
SK(I)=\{x\in I_K(G) \,\, \mid \,\, x^*=-x\}
$$
be  the  set of symmetric and  the  set of skew symmetric elements of
$I_K(G)$, respectively.

Clearly, $I_K(G)=S(I)+SK(I)$ and $S(I)\cap SK(I)=0$. Indeed, for any $y\in
I_K(G)$ we have $y=\frac{y+y^*}{2}+\frac{y-y^*}{2}$. 
Moreover, if $x\in S(I)$, then $x=\sum_{g\in G}\alpha_g(g+g\m1-2)$ and
$|S(I)|=|K|^{\frac{|G|-1}{2}}$, whence 
$$
|SK(I)|={|K|^{|G|-1}}:{|K|^\frac{{|G|-1}}{2}}=|K|^\frac{|G|-1}{2}.
$$ 

Recall that, if $k\in SK(I)$, then the element $1+k$ is a   unit
and the element $u=(1-k)(1+k)\m1$ is a   unitary unit \cite{6}, which is
called a  {\it Cayley unitary unit}.

It is easy to observe that  any element in $V_*(KG)$  is a 
Cayley  unitary unit for $char(K)>2$ . Indeed, if $u\in V_*(KG)$
then $1+u$ is a  unit because $\chi(1+u)=2\not=0$  and
$k=(1-u)(1+u)\m1$ is skew symmetric. Indeed, 
$$
k^*=(1+u^*)\m1(1-u)^*=(u\m1(1+u))\m1(u\m1(u-1))=-k.
$$
 Therefore, $1+k=((1+u)+(1-u))(1+u)\m1=2(1+u)\m1$ is
a unit  and $u=(1-k)(1+k)\m1$ is a Cayley unitary unit. We
conclude that $u=-1+2(1+k)\m1$ and the number of the unitary units 
equals  the number of skew symmetric elements and 
$$
|V_*(KG)|=|SK(I)|=|K|^\frac{|G|-1}{2}.  
$$

Now let us state  some basic properties of group algebras of the
$2$-groups. 
To determine the order of the unitary subgroup $V_*(KG)$,  
we need the following results from \cite{5}.

\proclaim
{Lemma 1}
Let $K$ be a field of prime characteristic $p$ and let $G$ be a
nonabelian locally finite $p$-group. The subgroup $V_*(KG)$ is
normal in $V(KG)$ if and only if $p=2$ and $G$ is the direct
product of an elementary abelian group with a group $H$ for which
one of the following holds\rom:
\itemitem{\rm(i)} $H$ has no direct factor of order $2$, but it
is a semidirect product of a group $\gp{h}$ of order $2$ and an
abelian $2$-group $A$  with $a^h=a\m1$ for all $a$ in $A$\rom;
\itemitem{\rm(ii)} $H$ is an extraspecial $2$-group  or the
central product of such a group with a cyclic group of order $4$.

\endproclaim
\rightline{\text{\qed}}

  Recall that a $p$-group is extraspecial if its centre,
commutator subgroup and Frattini subgroup  coincide  and have
order $p$.

\proclaim
{Lemma 2}
 Let $K$ be a commutative ring and $G$ be any group.  For $x\in
V(KG)$ and $y\in V_*(KG)$, we have $x\m1yx\in V_*(KG)$ if and
only if $xx^*$ commutes with $y$.
\endproclaim

\demo{Proof }
Clearly, $(x\m1yx)^*=(x\m1yx)\m1$ means that
$x^*y^*(x^*)\m1=x\m1y\m1x$  which in turn is equivalent to
$xx^*y^*=y\m1xx^*$.  Since we are given $y^*=y\m1$, this proves
the lemma.
\enddemo
\rightline{\text{\qed}}

Since  $G\leq V_*(KG)$, any  element,  which commutes with 
every element of $V_*(KG)$, is central in $KG$. Thus 
Lemma $2$  gives the following:

\proclaim
{Corollary 1}
 {The subgroup $V_*(KG)$ is normal in $V(KG)$ if and only if 
all elements of the form $xx^*$ with $x\in V(KG)$ are central 
in $KG$. }
\endproclaim

\rightline{\text{\qed}}

\proclaim
{Lemma 3}
 {Let $K$ be a field of characteristic $2$ and let $G$ be a nonabelian
locally finite $2$-group for which one of the following holds\rom:
\itemitem{\rm(i)} $G$  is a semidirect
product of a group $\gp{b}$ of order $2$ and an abelian $2$-group $A$, with
$a^b=a\m1$ for all $a$ in $A$\rom;
\itemitem{\rm(ii)} $G$ is an extraspecial $2$-group or the central product of
such a group with a cyclic group of order $4$.

\noindent
Then the map $x\to xx^*$ is a homomorphism of the group $V(KG)$ into the
subgroup $S_K(G)=\{xx^*\mid x\in V(KG)\}$ of the symmetric units and  if
$V(KG)$ is  finite, then the order
of the unitary subgroup $V_*(KG)$ coincides with the index of the subgroup
$S_K(G)$ in $V(KG)$.

}
\endproclaim

\demo{Proof}  
From Corollary 1 we have that $xx^*$ is central in $V(KG)$. Setting
$\phi(x)=xx^*$,  we have 
$$
\phi(xy)=(xy)(xy)^*=xyy^*x^*=(xx^*)(yy^*)=\phi(x)\phi(y) 
$$
for all $x,y\in V(KG)$. Therefore $\phi $ is an epimorphism and 
the kernel of $\phi$ is the  unitary subgroup $V_*(KG)$. From 
this follows the rest  of Lemma $3$.

\enddemo 
\rightline{\text{\qed}}

\subhead
3. The order of the unitary subgroup 
\endsubhead
For a finite $2$-group $G$ with the commutator  subgroup 
$G'=\gp{\, c \,  \mid\,  c^2=1\,}$  of order $2$ we define  
$L_G$ as  the subset of the elements of order $4$ such that 
$L_G\cap L_Gc$  is empty and $L_G\cup L_Gc$ coincide with
all elements of order $4$ of $G$. 

\proclaim
{Lemma 4}
Let $G$ be an extraspecial $2$-group of order $|G|=2^{2n+1}$ 
with $n\geq 2$. Then 
$$
|L_G|=2^{n-1}(2^n-(-1)^n).
$$
\endproclaim

\demo{Proof} 
If $G$ is an extraspecial $2$-group of order $|G|=2^{2n+1}$ 
with $n\geq 2$, then by Theorem 5.3.8 in \cite{7} we have 
$G=G_1\Y\cdots \Y G_n$, where $G_i$ is a quaternion group 
of order $8$. Then 
$\zeta(G)=G'=\gp{\,\,c \,\, \mid \,\, c^2=1\,}$  
and $G_i\cap G_j=\gp{c}$ for any $i\not=j$. Evidently every 
element of order $4$ of $G$ can be written as 
$$
x=z_{i_1}z_{i_2}\cdots z_{i_s},\tag1 
$$
where $z_{i_k}\in G_{i_k}$ has order $4$ and 
$i_1<i_2<\cdots <i_s$. Then the number $s$ is 
called the length of $x$. Clearly,  $z_{i_k}^2=c$  
and the length of the element of order $4$ is odd.

Let $H_k=H(i_1,i_2,\ldots, i_k)=G_{i_1}\Y G_{i_2}\Y \cdots \Y  G_{i_k}$, where $k$
is odd and  $i_1<i_2<\cdots <i_k$. We shall prove that there are 
precisely $3^k$ elements of   length $k$ in $L_{H_k}$. Of course, every  $L_{G_i}$
contains $3$ different elements and every element of length $k$ of the form (1)
has order $4$. We conclude that the number of elements of $L_{H_k}$ is
$3^k$. Since the number of different subgroups $H(i_1,i_2,\ldots, i_k)$ of $G$
is $\binom{n}{k}$, then  the number of  elements of $L_G$ equals
$
|L_G|=\sum_{j=1}^m\binom{n}{2j-1}3^{2j-1},
$
where $m=\frac{n+1}{2}$ if $n$ is odd and $m=\frac{n}{2}$ if $n$ is even.

According to the binomial theorem,  we have
$$
M_1=2^{2n}=(1+3)^n=\sum_{i=0}^{n}\binom{n}{i}3^{i},
$$
$$
M_2=(-1)^n2^{n}=(1-3)^n=\sum_{i=0}^{n}\binom{n}{i}(-1)^{i}3^{i},
$$
whence  
$$
|L_G|=\sum_{j=1}^{m }\binom{n}{2j-1}3^{2j-1}=\frac{1}{2}(M_1-M_2)=2^{n-1}(2^n-(-1)^n).
$$

\enddemo

\rightline{\text{\qed}}

\proclaim
{Lemma 5}
Let $G$ be a central product of an extraspecial $2$-group $H$ of order
$|H|=2^{2n+1}$ with  a cyclic group $\gp{d}$ of order $4$. Then   $|L_G|=2^{2n}$.
\endproclaim

\demo{Proof}
It is obvious that any element of order $4$ in $G$ either lies in $H$ or may
be written as $ud$, where $u\in H$ and $|u|\not=4$. The number of elements $u$ is
exactly $|H|-2|L_H|$. We conclude that 
$$
|L_G|=|L_H|+\frac{|H|-2|L_H|}{2}=\frac{|H|}{2}=2^{2n}.
$$

\enddemo
\rightline{\text{\qed}}

\proclaim
{Theorem 1}
Let $K$ be a finite field of characteristic $2$. 
\itemitem{(i)} If  $G$ is an  extraspecial $2$-group of order
$|G|=2^{2n+1}$ with $n\geq 2$, then 
$$
|V_*(KG)|=|K|^{2^{n-1}(2^{n+2}-2^n+(-1)^n)-1}.
$$ 
\itemitem{(ii)} If  $G$ is
 a central product of an extraspecial $2$-group $H$ of order
$|H|=2^{2n+1}$ with  a cyclic group $\gp{d}$ of order $4$, then
$$
|V_*(KG)|=|K|^{3\cdot 2^{2^n}-1}.
$$
\endproclaim

\demo{Proof} 
Recall that, if $x=\sum_{i=1}^t\alpha_ia_i\in V(KG)$, then 
$$
xx^*=1+\sum_{i,j \atop {i<j}}\alpha_i\alpha_j(a_ia_j\m1+a_ja_i\m1).
$$ 

If $a_ia_j\m1$ has order $2$, then $a_ia_j\m1=a_ja_i\m1$ and therefore the
support of $xx^*$ contain no elements of order $2$. Thus we obtain 
$$
xx^*=1+\sum_{b\in L_G}\alpha_bb(1+c)=\prod_{b\in L_G}(1+\alpha_bb(1+c)), \tag2
$$
where $\alpha_b\in K$. Of course, the number of elements of the form $xx^*$
with $x\in V(KG)$ is at most $|K|^{|L_G|}$. We shall prove  that the subgroup 
$S_K(G)=\{xx^* \,\mid \, x\in V(KG)\}$ has order $|K|^{|L_G|}$.

\rm{(i)} Now let $G$ be an extraspecial $2$-group of order
$|G|=2^{2n+1}$ with $n\geq 2$. Take  $b\in L_G$ and $\alpha_b\in K$. Then
there exists an  element $w_b\in G$ of order $2$, such that $(w_b,b)\not=1$. Indeed,
the length of $b=z_{i_1}z_{i_2}\cdots z_{i_{2k+1}}$ is odd. Then  we choose  
another $u_1\in G_{i_1}$ of order $4$ such that $(u_1, z_{i_1})\not=1$ and an
arbitrary element $u_2$ either from the set $\{z_{i_2}\cdots z_{i_{2k+1}}\}$
or, if $k=0$, from $G_t$ with $t\not=i_1$ and $|u_2|=4$. Then $w_b=u_1u_2$ is an
element of order $2$ and does not commute with $b$.

Since $(1+\alpha_b(b+w_b))(1+\alpha_b(b+w_b))^*=1+\alpha_bb(1+c)$, this shows us
that any factors of $\prod_{b\in L_G}(1+\alpha_bb(1+c))$ belong  to the
subgroup $S_K(G)$. Since $uu^*$ is central for arbitrary $u\in V(KG)$, we have
$$
\prod_{b\in L_G}(1+\alpha_b(b+w_b))(1+\alpha_b(b+w_b))^*=
(1+\alpha_{b_1}(b_1+w_{b_1}))\times
$$$$
\times
\big(
\prod_{b\in L_G\setminus \{b_1\}}(1+\alpha_b(b+w_b))(1+\alpha_b(b+w_b))^*
\big)
(1+\alpha_{b_1}(b_1+w_{b_1}))^*=
$$$$
=\prod_{b\in L_G}(1+\alpha_b(b+w_b))(\prod_{b\in L_G}(1+\alpha_b(b+w_b)))^*.
$$
Thus $|S_K(G)|=|K|^{|L_G|}$ and by  Lemmas $3$ and $4$ we
get  the result.

\rm{(ii)} Now let $G=H\Y \gp{d}$, where $H=G_1\Y\cdots \Y G_n$ is
an extraspecial $2$-group of order $|H|=2^{2n+1}$. If $|H|>8$ and
$b\in L_H$, then, as before, we can prove that
$1+\alpha_bb(1+c)\in S_K(G)$. Thus it remains to consider the
case, when the element $b$ of order $4$ commutes with any element
of order $2$ in $G$. In this case  there exists
$G_i=\gp{a_i,b_i}$ such that $\gp{G_i,b}=\gp{G_i}\Y \gp{b}$.
Evidently $x=1+a_i+b$ and $y=1+a_i+b_ib$ are units,
$xx^*=1+(a_i+b)(1+c)$ and $yy^*=1+a_i(1+c)$. Recall  that $xx^*$
and $yy^*$ are central units by Lemma $1$. Furthermore, it
follows that $$ 1+b(1+c)=xx^*yy^*=(xy)(xy)^*\in S_K(G).  $$ Thus,
we conclude that $|S_K(G)|=|K|^{|L_G|}$ and by Lemmas $3$ and $5$
we prove  (ii) of the theorem.

\enddemo
\rightline{\text{\qed}}

Let $A$ be an abelian group and we define $A[2]=\{a\in A \mid a^2=1\}$. 

\proclaim
{Theorem 2}
Let $K$ be a finite field of characteristic $2$ and $G$ has an abelian
subgroup $A$ of index $2$ and an  element $b$  which  inverts every 
element of $A$. 
\itemitem{(i)}  If  $b$  has order $2$, then 
$$
|V_*(KG)|=|K|^\frac{3|A|+|A[2]|-2}{2}.
$$
\itemitem{(ii)} If $b$  has order $4$, then 
$$
|V_*(KG)|=2\cdot|A^2[2]|\cdot |K|^{|A|+\frac{1}{2}|A[2]|-1}.
$$
\endproclaim

\demo{Proof} 
(i) Any element of $V(KG)$ can be written as $x=x_0+x_1b$, where $x_i\in KA$
    and $\chi(x_0)+\chi(x_1)=1$. Evidently, $x_i^b=x_i^*$ ($i=1,2$) and  
$$
xx^*=(x_0+x_1b)(x_0^*+x_1b)=x_0x_0^*+x_1x_1^*.
$$

Before observed that for any $y=\alpha_1a_1+\alpha_2a_2+\cdots+\alpha_sa_s\in
KA$ with $a_i\in A$ we have 
$$
yy^*=(\alpha_1^2+\alpha_2^2+\cdots+\alpha_s^2)+\sum_{i,j \atop{i<j}}\alpha_i\alpha_j(a_ia_j\m1+a_ja_i\m1).
$$

If $a_ia_j\m1$ is an element of order two, then $a_ia_j\m1=a_ja_i\m1$ and,
therefore, the support of the element $yy^*$ does not contain the element of order
$2$. As a consequence, $xx^*\in S_K(G)$ has a unique
expression in the form 
$$
1+\sum_{a\in L}\alpha_a(a+a\m1),\tag3
$$
where $L$ is a full system of representatives of the subset
$\{\, a\m1, a \}$ with $a\in A\setminus A[2]$.  Obviously, the
number of elements of $L$ is $l=\frac{|A|-|A[2]|}{2}$.  Hence the
order of the subgroup $S_K(G)$ is at most $|K|^l$. To prove that
it is the order of  $S_K(G)$,  it remains to note that for any
$z$ of the form (3) there exists $x\in V(KG)$ such that $xx^*=z$.
Indeed,  if
$$
z=1+\alpha_1(a_1+a_1\m1)+\alpha_2(a_2+a_2\m1)+\cdots+\alpha_s(a_s+a_s\m1)
$$
we put $x_0=\alpha_1a_1+\alpha_2a_2+\cdots+\alpha_sa_s$ and
$x_1=1+x_0$. Then $$
\chi(x_0)+\chi(x_1)=1+2(\alpha_1+\alpha_2+\cdots+\alpha_s)=1
$$
and $x=x_0+x_1b$ is a unit with  $xx^*=1+x_0+x_0^*=z$. Thus the
order of the subgroup $S_K(G)$ equals $|K|^l$. Using this result, it is easy
to find the order of the unitary subgroup , using  Lemma $3$.

(ii) Now let  $A$ be a finite abelian $2$-group. Then $V_*(KA)=A\times U$ by
\cite{3} and the group $V_*(KA)$ has order 
$$
|A^2[2]|\cdot |K|^{\frac{1}{2}(|A|+|A[2]|)-1}.
$$
 Let us  calculate the order of the
unitary subgroup for the group $G$ using the method of the paper \cite{3} .

Let $G$ has an abelian subgroup $A$ of
index $2$ and an  element $b$ of order $4$, such that $a^b=a\m1$ for all $a\in
A$. Consider the  subgroup 
$$
R=\{1+(1+b^2)zb \mid z\in KA\}=\prod_{u\in T}\gp{1+\lambda u(1+b^2)b\mid
\lambda\in K},
$$
where $T$ is a transversal to $\gp{b^2}$ in $A$. Clearly, the elements of $R$
are unitary units of order $2$ and $|K|^{\frac{|A|}{2}}$ is the order of $R$. 

Let $x\in V_*(KG)$.  Because $G\subseteq V_*(KG)$, either $x$ or $xb$ can be written as
$x=x_1(1+x_2b)$, where $x_i\in KA$ and $\chi(x_1)=1$. Then $x\in V_*(KG)$ if
and only if 
$$
x_1x_1^*(1+x_2x_2^*)=1
\,\,\,\,\,\,\,\,
\text{ and }
\,\,\,\,\,\,\,\,
x_2(1+b^2)=0.
$$ 
It follows that $x_2x_2^*=0$ and the element $x_1$ is unitary 
and $1+x_2b\in R$. Since 
$w\m1(1+\lambda z(1+b^2)b)w=1+\lambda w\m1w^* z(1+b^2)b\in R$, 
where $w\in V_*(KA)$  and
$R\cap U=\gp{1}$, we have that 
$W=U\ltimes  R$ is a subgroup of $V_*(KG)$. 
Finally, $b\m1wb=w^*=w\m1\in V_*(KA)$ for every $w\in V_*(KA)$ 
 and we conclude that
$V_*(KG)=G\ltimes (R \ltimes U)$ and 
$$
|V_*(KG)|=2\cdot|A^2[2]|\cdot |K|^{\frac{1}{2}(|A|+|A[2]|)-1}\cdot |K|^{\frac{|A|}{2}}=
2\cdot|A^2[2]|\cdot |K|^{|A|+\frac{1}{2}|A[2]|-1}.
$$
\enddemo
\rightline{\text{\qed}}

\proclaim
{Corollary 2}
Let $K$ be a finite field of characteristic $2$.  
\itemitem{(i)} If  
$D_{2^{n+1}}=\gp{\,\, a,b\,\, \mid \,\, a^{2^n}=b^2=1,
\,\,\, a^b=a\m1\,\, }$  is a dihedral group of order $2^{n+1}$, then
$$
|V_*(KD_{2^{n+1}})|=|K|^{3\cdot 2^{n-1}}.
$$
\itemitem{(ii)}  If  
$Q_{2^{n+1}}=\gp{\,\, a,b\,\, \mid \,\, a^{2^n}=1, a^{2^{n-1}}=b^2, \,\,\,
a^b=a\m1\,\, }$  
is a quaternion  group of order $2^{n+1}$, then
$$
|V_*(KQ_{2^{n+1}})|=4\cdot |K|^{2^n}.
$$ 
\endproclaim

\rightline{\text{\qed}}

\enddocument